\documentstyle[12pt]{article}
\newcommand{\Cos}{\mathop{\rm Cos}\nolimits}
\newcommand{\Sin}{\mathop{\rm Sin}\nolimits}

\newcommand{\kr}{\mathop{\rm ker}\nolimits}
\newcommand{\R}{\mathop{\rm Re}\nolimits}
\newcommand{\res}{\mathop{\rm res}\nolimits}
\newcommand{\sign}{\mathop{\rm sign}\nolimits}
\newcommand{\bg}{\phantom._0\Phi_1(-;0;q^2,i(1-q^2)q^2zs)}
\newcommand{\bgm}{\phantom._0\Phi_1(-;0;q^2,-i(1-q^2)q^2zs)}
\newcommand{\Ex}{E_{q^2}(-i(1-q^2)zs)}
\newcommand{\mEx}{E_{q^2}(i(1-q^2)zs)}

\newcommand{\eizso}{\ddag e_{q^2}(i(1-q^2)zs)\ddag }
\newcommand{\Eizso}{\ddag E_{q^2}(i(1-q^2)q^2zs)\ddag }

\newcommand{\Eizkm}{E_{q^2}(-i(1-q^2)z\xi)}

\newcommand{\p}{\partial}

\newtheorem{predl}{Proposition}[section]
\newtheorem{defi}{Definition}[section]

\newtheorem{cor}{Corollary}[section]
\newtheorem{lem}{Lemma}[section]

\newcommand{\beq}[1]{\begin{equation}\label{#1}}
\newcommand{\eq}{\end{equation}}
\def\mapright#1{\smash{
\mathop{\longrightarrow}\limits^{#1}}}
\def\mapleft#1{\smash{
\mathop{\longleftarrow}\limits^{#1}}}
\def\mapup#1{\Big\uparrow\rlap
{$\vcenter{\hbox{$\scriptstyle#1$}}$}}
\def\mapdown#1{\Big\downarrow\rlap
{$\vcenter{\hbox{$\scriptstyle#1$}}$}}

\begin{document}
\begin{flushright}

 ITEP-TH-68/97\\
\end{flushright}

\vspace{10mm}
\begin{center}
{\Large \bf The $q$-Fourier transform of $q$-distributions.}\\ 
\vspace{5mm}
M.A.Olshanetsky 
 \\ ITEP, 117259,
Moscow\\e-mail olshanet@heron.itep.ru\\
\vspace{5mm}
V.-B.K.Rogov
\\ MIIT, 101475, Moscow\\
e-mail m10106@sucemi.bitnet\\

\end{center}
\begin{abstract}
We consider functions on the lattice generated by the integer powers 
of $q^2$ for $0<q<1$ and construct the $q$-analog of Fourier transform
based on the Jackson integral in the space of distributions on the lattice.  
\end{abstract}
                   
\section{Introduction}
\setcounter{equation}{0}

The classical Fourier transform is a powerful tool  in the Harmonic analysis
 on the simple Lie groups and  
homogeneous spaces. In the quantum group case similar problems can be 
simplified by the Fourier transform in the quantum affine space. There exist
a few papers devoted to the quantum Fourier transform \cite{CZ,M},
 but their main object
is the quantum algebra of operators without specification of the spaces of
functions where these operators are defined. For this reason, applications of
 those results in concrete problems of Harmonic analysis are not so easy.
We investigate here 
the simplest one-dimensional case where the whole effect of the non-commutativity 
emerges only in
the lattice form of the space of function and a very simple algebra
 generated by functions
in the direct and inverse Fourier spaces. This particular situation can be 
considered as the $q$-deformation of the classical Fourier transform 
which attracted last time 
a particular attention due its relation to the generalizations of the classical
harmonic oscillator problem (see \cite{KS,AAS,RS,ARS}). While one of the primary 
interest of these papers is the investigations of different kernels of the 
$q$-Fourier transforms, we focus on  spaces
of functions and their special $q$-Fourier transform based on the $q$-exponents.  
Our approach is very similar to the classical picture developed in \cite{GS}. 
We consider the space of  test functions and $q^2$-distributions and then
 determine the 
$q^2$-Fourier transform of test functions.  A slightly different  $q$-Fourier transform on a subspace of test functions was
considered in \cite{K}. 
The $q^2$-Fourier transform 
of the $q^2$-distributions is determined by the Parseval identity.
   In conclusion, we present a table of the $q^2$-Fourier transforms of some 
$q^2$-distributions.

\section{Preliminary relations}
\setcounter{equation}{0}

We assume that $z\in{\bf C}$ and $|q|<1$, unless otherwise is specified.

We recall some notations \cite{GR}. For an arbitrary $a$ 
$$ 
(a,q)_n=\left\{ 
\begin{array}{lcl} 
1&for&n=0\\
(1-a)(1-aq)\ldots(1-aq^{n-1})&for&n\ge1,\\
\end{array}
\right.
$$
$$
(a,q)_\infty=\lim_{n\to\infty}(a,q)_n.
$$
$$
\left[\begin{array}{c}
l \\i
\end{array}
\right ]_{q^2}
=\frac
{(q^2,q^2)_l}
{(q^2,q^2)_i(q^2,q^2)_{l-i}}.
$$
Consider the $q^2$-exponentials
\begin{equation}\label{2.3}
e_{q^2}(z)=\sum_{n=0}^\infty\frac{z^n}{(q^2,q^2)_n}=
\frac1{(z,q^2)_\infty},\qquad|z|<1,
\end{equation}
\begin{equation}\label{2.4}
E_{q^2}(z)=\sum_{n=0}^\infty\frac{q^{n(n-1)}z^n}{(q^2,q^2)_n}=
(-z,q^2)_\infty
\end{equation}
and the $q^2$-trigonometric functions
\begin{equation}\label{2.5}
\cos_{q^2}z=\frac12(e_{q^2}(iz)+e_{q^2}(-iz)),\quad
\sin_{q^2}z=\frac1{2i}(e_{q^2}(iz)-e_{q^2}(-iz)),
\end{equation}
\begin{equation}\label{2.6}
\Cos_{q^2}z=\frac12(E_{q^2}(iz)+E_{q^2}(-iz)),\quad
\Sin_{q^2}z=\frac1{2i}(E_{q^2}(iz)-E_{q^2}(-iz)).
\end{equation}
We also consider the basic hypergeometric series
\begin{equation}\label{2.7}
\phantom._0\Phi_1(-;0;q^2,z)=
\sum_{n=0}^\infty\frac{q^{2n(n-1)}z^n}{(q^2,q^2)_n}.
\end{equation}
\bigskip

It follows from (\ref{2.3}) that $e_{q^2}(z)$ is a meromorphic
function with simple poles in the points $z=q^{-2k}, 
k=0,1,\ldots.$ 
\begin{lem}\label{l2.1} 
The function $e_{q^2}(z)$ is represented as the sum of partial 
fractions:  
\begin{equation}\label{2.8}
e_{q^2}(z)=\frac1{(q^2,q^2)_\infty}\sum_{k=0}^\infty
\frac{(-1)^kq^{k(k+1)}}{(q^2,q^2)_k(1-zq^{2k})}.
\end{equation}
\end{lem}
{\bf Proof.} Let
$$
e_{q^2}(z,n)=\frac1{(z,q^2)_n}=\sum_{k=0}^n\frac{c_{k,n}}
{1-zq^{2k}},
$$
where
$$
c_{k,n}=\res_{z=q^{-2k}}e_{q^2}(z,n)=\lim_{z\to q^{-2k}}
(1-zq^{2k})e_{q^2}(z,n)=
$$
$$
=\frac1{(1-q^{-2k})\ldots(1-q^{-2})(1-q^2)\ldots(1-q^{2(n-k)})}=
\frac{(-1)^kq^{k(k+1)}}{(q^2,q^2)_k(q^2,q^2)_{n-k}}.
$$
Therefore,
$$
e_{q^2}(z,n)=\sum_{k=0}^\infty\frac{(-1)^kq^{k(k+1)}}
{(q^2,q^2)_k(q^2,q^2)_{n-k}(1-aq^{2k})}.
$$
From $\frac1{(q^2,q^2)_{n-k}}<\frac1{(q^2,q^2)_\infty}=
e_{q^2}(q^2)$, we obtain
$$
e_{q^2}(z)=\lim_{n\to\infty}e_{q^2}(z,n)=\frac1{(q^2,q^2)_\infty}
\sum_{k=0}^\infty\frac{(-1)^kq^{k(k+1)}}{(q^2,q^2)_k(1-zq^{2k})}.
$$
Since this series converges absolutely for $q\in(0,1)$, 
(\ref{2.8}) holds for an arbitrary $z\ne q^{-2k}$.  
\rule{5pt}{5pt}

From (\ref{2.5}) we obtain
\begin{cor}\label{c2.1}
The functions $\cos_{q^2}z$ and $\sin_{q^2}z$ are represented as the sums
of partial fractions
\begin{equation}\label{2.10}
\cos_{q^2}z=\frac1{(q^2,q^2)_\infty}\sum_{k=0}^\infty
\frac{(-1)^kq^{k(k+1)}}{(q^2,q^2)_k(1+z^2q^{4k})},
\end{equation}
\begin{equation}\label{2.11}
\sin_{q^2}z=\frac z{(q^2,q^2)_\infty}\sum_{k=0}^\infty
\frac{(-1)^kq^{k(k+3)}}{(q^2,q^2)_k(1+z^2q^{4k})}.
\end{equation}
\end{cor}
\begin{predl}\label{p2.1}
If $z$ and $q$ are real then  
$$ 
|\cos_{q^2}z|\le\frac{(-q^2,q^2)_\infty}{(1+z^2)(q^2,q^2)_\infty},
$$
$$
|\sin_{q^2}z|\le\frac{|z|(-1,q^2)_\infty}{(1+z^2)(q^2,q^2)_\infty}.
$$
\end{predl}
{\bf Proof.} Since
$$
\frac1{1+z^2q^{4k}}\le\frac{q^{-4k}}{1+z^2}, 
$$
for an arbitrary real $z$ and for an arbitrary $k\ge0$, we have
$$
|\sum_{k=0}^\infty\frac{(-1)^kq^{k(k+1)}}{(q^2,q^2)_k(1+z^2q^{4k})}|\le
\frac1{1+z^2}\sum_{k=0}^\infty\frac{q^{k(k-1)}q^{-2k}}
{(q^2,q^2)_k}=\frac1{1+z^2}E_{q^2}(q^{-2}),
$$
$$
|\sum_{k=0}^\infty\frac{(-1)^kq^{k(k+3)}}{(q^2,q^2)_k(1+z^2q^{4k})}|\le
\frac1{1+z^2}\sum_{k=0}^\infty\frac{q^{k(k-1)}}{(q^2,q^2)_k}=
\frac1{1+z^2}E_{q^2}(1).
$$
The statement of Proposition follows from (\ref{2.10}), 
(\ref{2.11}) and the last relations.  \rule{5pt}{5pt} 

The next proposition follows immediately from  (\ref{2.4}) and 
(\ref{2.6}).  
\begin{predl}\label{p2.2} $\Cos_{q^2}z$ and $\Sin_{q^2}z$ 
satisfy the inequalities 
\begin{equation}\label{2.14} 
|\Cos_{q^2}z|\le1,\qquad |\Sin_{q^2}z|\le|z|.  
\end{equation} 
\end{predl} 

Consider the theta function 
$$
\vartheta(u|\tau)=2p^{1/8}\sin\pi u
\prod_{n>0}(1-e^{2\pi i u}p^n)(1-e^{-2\pi i u}p^n)(1-p^n)
$$
where  $p=\exp 2\pi i \tau,~~(Im\tau>0)$. 
Let 
\begin{equation}\label{2.15}
{\bf Q}(z,q)=(1-q^2)\sum_{m=-\infty}^\infty\frac1
{zq^{2m}+z^{-1}q^{-2m}}
\end{equation}
Then it follows from the modular properties of (\ref{2.15}) that
$$
{\bf Q}(z,q)=-\frac{1}{2\pi i}
\frac
{\vartheta(u+\tau/2|\tau)\vartheta'(0|\tau)}
{\vartheta(u|\tau)\vartheta(\tau/2|\tau)},
$$
where 
$$
\tau=-\frac{2i\ln q}{\pi},~~u=\frac{1}{\pi i}\ln z+\frac{1}{2}.
$$

Define the function
\begin{equation}\label{2.16}
\Theta(z)=(1-q^2)\sum_{m=-\infty}^\infty\sin_{q^2}((1-q^2)q^{2m}z).
\end{equation}
\begin{predl}\label{p2.3}
 $\Theta(z)$ has the following properties:\\
1. $\Theta(q^{2k}z)=\Theta(z)$ for any $z\ne0$,\\
2. $\Theta(z)={\bf Q}((1-q^2)z,q).$
\end{predl}
{\bf Proof.}  1. follows immediately from (\ref{2.16}).

Due to (\ref{2.11}) $\Theta(z)$ can be represented as
$$
\Theta(z)=\frac{1-q^2}{(q^2,q^2)_\infty}\sum_{m=-\infty}^\infty
q^{2m}\sum_{k=0}^\infty\frac{(-1)^kq^{k(k+3)}(1-q^2)z}
{(q^2,q ^2)_k[1+(1-q^2)^2z^2q^{4(m+k)}]}.
$$
We can
change the order of summation 
since the inner series converges uniformly with respect to $m$. 
Then using (\ref{2.3}), (\ref{2.4}) 
and (\ref{2.15}) we obtain 
$$ 
\Theta(z)=\frac{1-q^2}{(q^2,q^2)_\infty}\sum_{k=0}^\infty 
\frac{(-1)^kq^{k(k+3)}}{(q^2,q^2)_k}\sum_{m=-\infty}^\infty
\frac{(1-q^2)zq^{2m}}{1+(1-q^2)^2z^2q^{4(m+k)}}=
$$
$$
=\frac{1-q^2}{(q^2,q^2)_\infty}\sum_{k=0}^\infty
\frac{(-1)^kq^{k(k+1)}}{(q^2,q^2)_k}\sum_{m=-\infty}^\infty
\frac{(1-q^2)zq^{2m}}{1+(1-q^2)^2z^2q^{4m}}={\bf Q}((1-q^2)z,q).
$$
\rule{5pt}{5pt}

Specially we define
\begin{equation}\label{2.17} 
\Theta_0=\Theta(1)={\bf Q}(1-q^2,q). 
\end{equation} 
\begin{predl}\label{p2.4} 
For an arbitrary integer $M>0$ 
\begin{equation}\label{2.18} 
(1-q^2)z\sum_{m=-M}^\infty 
q^{2m}\cos_{q^2}((1-q^2)q^{2m}z)= \sin_{q^2}((1-q^2)q^{-2M}z), 
\end{equation}
\begin{equation}\label{2.19}
(1-q^2)z\sum_{m=-M}^\infty q^{2m}\sin_{q^2}((1-q^2)q^{2m}z)=
1-\cos_{q^2}((1-q^2)q^{-2M}z).
\end{equation}
\end{predl}
{\bf Proof.} Assume  that $|(1-q^2)q^{-2m}z|<1$. Then
$$
(1-q^2)z\sum_{m=-M}^\infty q^{2m}\cos_{q^2}((1-q^2)q^{2m}z)=
\sum_{m=-M}^\infty(1-q^2)zq^{2m}\sum_{k=0}^\infty(-1)^k
\frac{(1-q^2)^{2k}q^{4mk}z^{2k}}{(q^2,q^2)_{2k}}
$$
\begin{equation}\label{2.20}
=\sum_{k=0}^\infty(-1)^k\frac{(1-q^2)^{2k+1}z^{2k+1}}
{(q^2,q^2)_{2k}}\sum_{m=-M}^\infty q^{2m(2k+1)}=
\sin_{q^2}((1-q^2)q^{-2M}z).
\end{equation}
For $z\ne\pm i(1-q^2)^{-1}q^{-2n}$ the left side of (\ref{2.18}) 
is defined as the analytic continuation of (\ref{2.20}).  
(\ref{2.19}) can be proved in the same way. \rule{5pt}{5pt}

\section{$q^2$-derivatives and $q^2$-integral}
\setcounter{equation}{0}
Let ${\cal A}=C(z,z^{-1})$ be the algebra of formal Laurent series.  
The $q^2$-derivative of $f(z)\in {\cal A}$ is defined as
\begin{equation}\label{3.1} 
\p_zf(z)=\frac{z^{-1}}{1-q^2}(f(z)-f(q^2z)).  
\end{equation} 
For an arbitrary $n\ge0$
\begin{equation}\label{3.2}
\p_z^kz^n=\left\{
\begin{array}{lcl}
\frac{(q^2,q^2)_n}{(q^2,q^2)_{n-k}(1-q^2)^k}z^{n-k}&for&0\le k\le n\\
0&for&k>n,\\
\end{array}
\right.
\end{equation}
and for any $n\ge0$ and $k\ge0$
\begin{equation}\label{3.3}
\p_z^kz^{-n-1}=(-1)^kq^{-k(2n+k+1)}\frac{(q^2,q^2)_{n+k}}
{(q^2,q^2)_n(1-q^2)^k}z^{-n-k-1}.
\end{equation}

The $q^2$-integral (Jackson integral \cite{GR})  
is defined as the map $I_{q^2}$ from ${\cal A}$ to the space 
of formal number series
\begin{equation}\label{3.4} 
I_{q^2}f=\int d_{q^2}zf(z)=(1-q^2)\sum_{m=-\infty}^\infty q^{2m} 
[f(q^{2m})+f(-q^{2m})]  
\end{equation}
\begin{defi}\label{d3.2}
The function $f(z)$ is locally $q^2$-integrable, 
if the $q^2$-integral 
\begin{equation}\label{3.5} 
\int_a^bd_{q^2}zf(z)=(1-q^2)\sum_{m=0}^\infty q^{2m}
[bf(bq^{2m})-af(aq^{2m})]
\end{equation}
exists for any finite $a$ and $b$, i.e., the series in the 
right side of (\ref{3.5}) converges.  
\end{defi}            

\begin{defi}\label{d4.4}
 $f(z)$ is absolutely $q^2$-integrable function, if the series
$$
\sum_{m=-\infty}^\infty q^{2m}[|f(q^{2m})|+|f(-q^{2m})|]
$$
converges.
\end{defi}

Let ${\cal B}$ be the algebra analogous to  ${\cal A}$, 
but generated by $s,s^{-1}$ which commute with $z$ 
as $zs=q^2sz$. Define the 
$q^2$-differentiation in ${\cal B}$ as 
$$ 
\p_s\phi(s)=(\phi(s)-\phi(q^2s))\frac{s^{-1}}{1-q^2}.
$$ 

We denote by ${\cal AB}$ the whole algebra with two generators 
$z,z^{-1},s,s^{-1}$ and the \\ $q^2$-differentiation
\begin{equation}\label{3.7}
zs=q^2sz,\quad \p_zs=q^{-2}s\p_z,\quad \p_sz=q^2z\p_s,
\quad \p_z\p_s=q^2\p_s\p_z.
\end{equation}
 
We consider ${\cal AB}$ as a left module under the left action of
${\cal A}$ by multiplication, and a right module 
under the left action of ${\cal B}$.

To define the $q^2$-integral on ${\cal AB}$ we order the generators
of integrand
in such a way that $z$ stays on the left side while $s$ 
stays on the right side.
For example, if $f(z)=\sum_ra_rz^r$ then
$$
f(zs)=\sum_ra_r(zs)^r=\sum_ra_rq^{-r(r-1)}z^rs^r.
$$
For convenience we introduce the following notation
$$
\ddag g(zs)\ddag=\sum_ra_rz^rs^r,~{\rm if}~g(z)=\sum_ra_rz^r.
$$
 For example, we can derive from 
(\ref{2.3}), (\ref{2.4}) and (\ref{3.7}) that
$$ 
\mEx=\eizso.  
$$
Following prescribed rules we calculate
\\ 
1.  
\begin{equation}\label{3.8}
\int d_{q^2}zz^{-1}\mEx=2i\Theta(s)
\end{equation}
(see (\ref{2.17})).
\\
2. Using (\ref{2.19}) and Proposition \ref{p2.1} we get
$$
\int d_{q^2}z\Theta(z)\mEx=
2i\Theta_0s^{-1}.
$$

\section{ $q^2$-distributions}
\setcounter{equation}{0}

Let $S_{q^2}=\{\phi(x)\}$ be the space of infinitely $q^2$-differentiable 
fast decreasing functions 
\begin{equation}\label{4.1} 
|x^k\p_x^l\phi(x)|\le C_{k,l}(q),~k\ge0,~ l\ge0. 
\end{equation}

Let $S$ be the space of infinitely differentiable (in the classical sense)
fast decreasing functions  
\begin{equation}\label{4.2} 
|x^k\phi^{(l)}(x)|\le C_{k,l},~k\ge0~l\ge0.  
\end{equation}
\begin{predl}\label{p4.1}
$$
S\subset S_{q^2}.
$$
\end{predl}
{\bf Proof.} We denote by $\phi(x_1,x_2,\ldots,x_k)$ 
the separated difference of order $k$. Then it is easy to find 
that for an arbitrary integer $l\ge0$ 
\begin{equation}\label{4.3} 
\phi(q^{2l}x,q^{2l-2}x,\ldots,x)=
\frac{(q^2,q^2)_l}{(1-q^2)^l}\p_x^l\phi(x).
\end{equation}
On the other hand, if $\phi(x)$ is differentiable $l$ times (in the 
classical sense), then there exists a point $\xi\in(q^{2l}x,x)$,
 such that 
\begin{equation}\label{4.4} 
\phi(q^{2l}x,q^{2l-2}x,\ldots,x)=\frac1{l!}\phi^{(l)}(\xi).
\end{equation}
Then, it follows from (\ref{4.3}) and (\ref{4.4})
$$
\p_x^l\phi(x)=\frac{(1-q^2)^l}{(q^2,q^2)_ll!}\phi^{(l)}(\xi),\quad
\xi\in(q^{2l}x,x).
$$
Therefore, from (\ref{4.2})
$$
|x^k\p_x^l\phi(x)|\le\frac{(1-q^2)^l}{(q^2,q^2)_ll!}C_{k,l}.
$$
Finally, if $\phi_n(x)\to0$ for $n\to\infty$ in the space $S$ then 
$\phi_n(x)\to0$ in the space $S_{q^2}$. \rule{5pt}{5pt}

\begin{defi}\label{d4.1} 
The skeleton $\hat\phi(z)$ of $\phi(z)\in S_{q^2}$ is the set of
evaluations of $\phi(z)$ on the lattice generated by the
powers of $q^2$ 
$$ 
\hat\phi(z)=\phi(z)|_{z=q^{2n}}, n=0,\pm1,\ldots.  
$$ 
\end{defi}

The space of the skeletons we denote by $\hat S_{q^2}$.

The functions
$$
\hat\phi_n^+(z)=\left\{\begin{array}{rcl}
1 & {\rm for}& z=q^{2n}\\ 0 & {\rm for} & z\ne q^{2n},\\
\end{array}
\right.
\hat\phi_n^-(z)=\left\{\begin{array}{rcl}
1&{\rm for}&z=-q^{2n}\\ 0&{\rm for}&z\ne-q^{2n}.\\
\end{array}
\right.
$$
generate the basis in the space $\hat S_{q^2}$:
$$
\hat\phi(z)=\sum_{n=-\infty}^\infty[a_n\hat\phi_n^+(z)+
b_n\hat\phi_n^-(z)],
\quad a_n=\phi(q^{2n}),\quad b_n=\phi(-q^{2n}).
$$

The topology in $\hat S_{q^2}$ is induced by the topology in 
$ S_{q^2}$ (\ref{4.1}):
$$
|\frac{q^{2n(k-l)}}{(1-q^2)^l}\sum_{i=1}^l(-1)^i
\left[\begin{array}{c}
l \\i
\end{array}
\right ]_{q^2}
q^{-i(2l-i-1)}a_{n+i}|\le C_{k,l}(q),
$$
$$
|(-1)^{k-l}\frac{q^{2n(k-l)}}{(1-q^2)^l}\sum_{i=1}^l(-1)^i
\left[\begin{array}{c}
l \\i
\end{array}
\right ]_{q^2}
q^{-i(2l-i-1)}b_{n+i}|\le C_{k,l}(q).
$$

Let $L_{q^2}$ be the linear map from $S_{q^2}$ to $\hat S_{q^2}$ 
defined
by the evaluation of functions in the vertices of the lattice. Then
$$  
\hat S_{q^2}=S_{q^2}/\kr L_{q^2}.   
$$

Let $\Lambda\phi(z)=\phi(q^2z)$.
Then
\begin{equation}\label{4.6}
\Lambda z=q^2z\Lambda,\qquad \p_z\Lambda=q^2\Lambda\p_z. 
\end{equation}
The operations $\Lambda$ and
$\p_z$  are well defined
on $\hat S_{q^2}$ since
$$
\Lambda L_{q^2}=L_{q^2}\Lambda,~~~
\p_zL_{q^2}=L_{q^2}\p_z.
$$
Moreover, the $q^2$-integral vanishes on the functions from
$\kr L_{q^2}$, and therefore is defined on the quotient
$
\hat S_{q^2}=S_{q^2}/\kr L_{q^2}.
$
In some places, we will write in the $q^2$-integral an element from 
$\hat S_{q^2}$ assuming that it is a representative from the quotient.
\begin{predl}\label{p4.4} 
If $\hat\phi(z)\in\hat S_{q^2}$ then 
$$
\int d_{q^2}z\p_z\hat\phi(z)=0.
$$ 
\end{predl} 
{\bf Proof.} It follows from (\ref{3.1}),(\ref{3.4}) and (\ref{4.1}) 
$$ 
\int d_{q^2}z\p_z\hat\phi(z)=\sum_{m=-\infty}^\infty
[\hat\phi(q^{2m})- 
\hat\phi(q^{2m+2})-\hat\phi(-q^{2m})+\hat\phi(-q^{2m+2})]=
$$
$$
=\lim_{M\to\infty}[\hat\phi(q^{-2M})-\hat\phi(-q^{-2M})]=0.
$$ 
\rule{5pt}{5pt}
\begin{cor}\label{c4.1}
($q^2$-integration by parts) For any $k\ge 0$ 
$$ 
\int d_{q^2}z\hat\phi(z)\p_z^k\hat\psi(z)=
(-1)^kq^{-k(k-1)} \int d_{q^2}z\p_z^k\hat\phi(z)\hat\psi(q^{2k}z).  
$$ 
\end{cor}

\begin{defi}\label{d4.3}
The $q^2$-distribution $f$ over $\hat S_{q^2}$ 
is a linear continuous functional $f:\hat S_{q^2}\to{\bf C}.$ 
\end{defi}

We denote by $\hat S_{q^2}'$ the space of the $q^2$-distributions 
over $\hat S_{q^2}$.

\begin{defi}\label{d4.5}
The sequence $\{f_n\}$  is referred to as convergent in  
$\hat S_{q^2}'$ to $f\in\hat S_{q^2}'$, if for any $\phi(z)\in\hat 
S_{q^2}$ the  sequence $\{<f_n,\phi>\}$ converges to $<f,\phi>$.  
\end{defi} 

The $q^2$-distributions which are defined by the $q^2$-integral
$$
<f,\phi>=\int_\infty^\infty d_{q^2}z\bar f(z)\phi(z)=
(1-q^2)\sum_{m=-\infty}^\infty q^{2m} 
[\bar f(q^{2m})\phi(q^{2m})+\bar f(-q^{2m})\phi(-q^{2m})] 
$$ 
we refer as regular.

 Corollary \ref{c4.1} and Proposition \ref{p4.4} allows to
introduce the $q^2$-differentiation in $\hat S_{q^2}'$
\begin{equation}\label{4.9}
<\p_zf,\phi>=- q^2<\Lambda f,\p_z\phi>.
\end{equation}

Let us give some examples of $q^2$-distributions.\\
1.
\begin{equation}\label{4.10}
<\theta_{q^2}^+,\phi>=\int_0^\infty d_{q^2}z\hat\phi(z)=
(1-q^2)\sum_{m=-\infty}^\infty q^{2m}\phi(q^{2m}).
\end{equation}
2.
\begin{equation}\label{4.11}
<\theta_{q^2}^-,\phi>=\int_{-\infty}^0d_{q^2}z\hat\phi(z)=
(1-q^2)\sum_{m=-\infty}^\infty q^{2m}\phi(-q^{2m}).
\end{equation}
So, the $q^2$-distributions $\theta_{q^2}^+$ and $\theta_{q^2}^-$ 
correspond to the functions
$$
\hat\theta^+(z)=\sum_{n=-\infty}^\infty\hat\phi_n^+(z),\qquad
\hat\theta^-(z)=\sum_{n=-\infty}^\infty\hat\phi_n^-(z).
$$
\\ 
3.  
$$ 
<\delta_{q^2},\phi>=\phi(0)=\lim_{m\to\infty}
\frac{\phi(q^{2m})+\phi(-q^{2m})}2.
$$
Moreover, it follows from (\ref{4.9})
$$
<\p_z(\theta_{q^2}^+-\theta_{q^2}^-),\phi>=
-<\theta_{q^2}^+-\theta_{q^2}^-,\p_z\phi>=
$$
$$
=-(1-q^2)\sum_{m=-\infty}^\infty q^{2m}\Bigl(\frac{q^{-2m}}{1-q^2}
[\phi(q^{2m})-\phi(q^{2m+2})]+\frac{q^{-2m}}{1-q^2}
[\phi(-q^{2m})-\phi(-q^{2m+2})]\Bigr)=
$$
$$
=-\lim_{M\to\infty}\sum_{m=-M}^M[\phi(q^{2m})-\phi(q^{2m+2})+
\phi(-q^{2m})-\phi(-q^{2m+2})]=
$$
$$
=-\lim_{M\to\infty}[\phi(q^{-2M})-\phi(q^{2M+2})+
\phi(-q^{-2M})-\phi(-q^{2M+2})]=2\phi(0),
$$
i.e.,
\begin{equation}\label{4.13}
\p_z(\theta_{q^2}^+(z)-\theta_{q^2}^-(z))=2\delta_{q^2}(z).
\end{equation}
4.
$$
<z^{-1},\phi>=\int_0^\infty d_{q^2}zz^{-1}[\hat\phi(z)-\hat\phi(-z)]=
(1-q^2)\sum_{m=-\infty}^\infty[\phi(q^{2m})-\phi(-q^{2m})].
$$
5. For an arbitrary $k\ge0$
$$
<z^{-k-1},\phi>=(-1)^kq^{k(k+1)}\frac{(1-q^2)^k}{(q^2,q^2)_k}
<\p_z^kz^{-1},\phi>=\frac{(1-q^2)^k}{(q^2,q^2)_k}<z^{-1},\p_z^k\phi>=
$$
\begin{equation}\label{4.15}
=\frac{(1-q^2)^{k+1}}{(q^2,q^2)_k}\sum_{m=-\infty}^\infty
[\p_z^k\phi(z)|_{z=q^{2m}}-\p_z^k\phi(z)|_{z=-q^{2m}}].
\end{equation}
6. For an arbitrary $\nu>-1$
$$
<z_+^\nu,\phi>=\int_0^\infty d_{q^2}zz^\nu\hat\phi(z)=
(1-q^2)\sum_{m=-\infty}^\infty q^{2m(\nu+1)}\phi(q^{2m}).
$$
Since for any $k\ge0$
$$
\p_z^kz^\nu=(-1)^kq^{k(2\nu-k+1)}\frac{(q^{-2\nu},q^2)_k}{(1-q^2)^k}
z^{\nu-k},
$$
we define 
$$
<z_+^{\nu-k},\phi>=\frac{(1-q^2)^{k+1}}{(q^{-2\nu},q^2)_k}
\sum_{m=-\infty}^\infty q^{2m(\nu+1)}\p_z^k\phi(z)|_{z=q^{2m}}.
$$
7. Similarly, for an arbitrary $\nu>-1$
$$
<z_-^\nu,\phi>=\int_{-\infty}^0d_{q^2}z(-z)^\nu\hat\phi(z)=
(1-q^2)\sum_{m=-\infty}^\infty q^{2m(\nu+1)}\phi(-q^{2m}),
$$
and for an arbitrary $k\ge0$
$$
<z_-^{\nu-k},\phi>=(-1)^k\frac{(1-q^2)^{k+1}}{(q^{-2\nu},q^2)_k}
\sum_{m=-\infty}^\infty q^{2m(\nu+1)}\p_z^k\phi(z)|_{z=-q^{2m}}.
$$

\section{The $q^2$-Fourier transform} 
\setcounter{equation}{0}

Consider  the space $ S^{q^2}=\{\psi(s)\}$ generated 
by the similar class  of functions  as $ S_{q^2}=\{\phi(z)\}$ (\ref{4.1}), 
but $s$ and $z$ are related as the generators in the algebra ${\cal AB}$
(\ref{3.7}). We introduce the same topology in $S^{q^2}$ as in 
$S_{q^2}$
\begin{equation}\label{top}
|s^k\p_s^l\phi(s)|\le C_{k,l}(q),~k\ge0,~ l\ge0 ,
\end{equation}
and, thereby, these spaces are isomorphic.

We define the map from $S_{q^2}$ to  $S^{q^2}$
$$
\def\normalbaselines{\baselineskip20pt
\lineskip3pt     \lineskiplimit3pt}
\matrix{S_{q^2}  &
\mapright{L_{q^2}} &
 \hat S_{q^2} &
\mapright{{\cal F}_{q^2}}&
   S^{q^2},\cr}
$$
where
\begin{equation}\label{5.1}
{\cal F}_{q^2}\phi(z)=\int d_{q^2}z\phi(z)\bg,
\end{equation}
is the $q^2$-Fourier transform and 
$\phantom._0\Phi_1$ is defined by (\ref{2.7}). 
It follows from (\ref{5.1}) that ${\cal F}_{q^2}$ is well-defined
on $S_{q^2}$, and, in some places, we will preserve the same notation 
for the operator acting from $S_{q^2}$ to $S^{q^2}$ 
discarding the action of $L_{q^2}$.

Our aim is to find the inverse map and to prove the continuity of the
both maps. More precisely, we define the dual
space of skeletons $\hat{S}^{q^2}$   for $S^{q^2}$ and the map
\begin{equation}\label{5.4} 
{\cal F}_{q^2}^{-1}\psi(s)=\frac1{2\Theta_0}\int\Ex\psi(s)d_{q^2}s,
~~\psi(s)\in\hat{S}^{q^2}
\end{equation}
from $\hat{S}^{q^2}$ to $S_{q^2}$. 

Consider the diagram
$$
\def\normalbaselines{\baselineskip20pt
\lineskip3pt     \lineskiplimit3pt}
\matrix{
S_{q^2}                  & \mapright{L_{q^2}} & \hat S_{q^2} \cr
\mapup{{\cal F}_{q^2}^{-1}}&                       
&\mapdown{{\cal F}_{q^2}}\cr
\hat S^{q^2}                        & \mapleft{L_{q^2}} & 
       S^{q^2}. \cr
}
$$
We will prove
\begin{predl}\label{p5.1}
1.The maps 
$$
{\cal F}_{q^2}L_{q^2} : S_{q^2}\rightarrow S^{q^2}
$$
$$
{\cal F}^{-1}_{q^2}L_{q^2} : S^{q^2}\rightarrow S_{q^2}
$$
realize isomorphisms of topological spaces.\\
2.The maps
$$
L_{q^2}{\cal F}_{q^2}^{-1}L_{q^2}{\cal F}_{q^2}:
\hat S_{q^2}\rightarrow\hat S_{q^2},
$$
$$
L_{q^2}{\cal F}_{q^2}L_{q^2}{\cal F}_{q^2}^{-1}
$$
are the identity maps on $\hat S_{q^2}$ and $\hat S^{q^2}$ correspondingly.
\end{predl}

We begin with
\begin{lem}\label{l5.1}
$$
\int d_{q^2}ze_{q^2}(-i(1-q^2)z)\bg=
\left\{
\begin{array}{lcl}
\frac2{1-q^2}\Theta_0 & for & s=1\\ 0 & for & s\ne 1.\\
\end{array}
\right.
$$
where $\Theta_0$ is determined by (\ref{2.17}).
\end{lem}
{\bf Proof.} In according 
with our definition of  $q^2$-integrals all integrands must by ordered.
It follows from (\ref{2.4}), (\ref{2.7}) and (\ref{3.7}) 
 that 
$$ 
\bg=\Eizso,  
$$ 
and we come to the integral
$$ 
\int d_{q^2}ze_{q^2}(-i(1-q^2)z)\Eizso.
$$ 
 Using (\ref{2.3}) and 
(\ref{2.4}) we obtain
$$ 
\p_z[e_{q^2}(-i(1-q^2)z)\ddag\mEx\ddag]= -ie_{q^2}(-i(1-q^2)z)\Eizso (1-s).  
$$ 
Hence, if $s\ne1$, then 
$$ 
\int d_{q^2}ze_{q^2}(-i(1-q^2)z)\Eizso= 
$$ 
$$ 
=i(1-s)^{-1}\int d_{q^2}z\p_z[e_{q^2}(-i(1-q^2)z)\ddag\mEx\ddag].  
$$ 
 It can be represented 
in the form (see(\ref{2.5}) and (\ref{2.6}))
$$ 
(s-1)^{-1}\lim_{M\to\infty} 
[\cos_{q^2}((1-q^2)q^{-2M})\Sin_{q^2}((1-q^2)q^{-2M+2}s)-
$$
$$
-\sin_{q^2}((1-q^2)q^{-2M})\Cos_{q^2}((1-q^2)q^{-2M+2}s)-
\cos_{q^2}((1-q^2)q^{2M+2})\Sin_{q^2}((1-q^2)q^{2M+4}s)+
$$
$$
+\sin_{q^2}((1-q^2)q^{2M+2})\Cos_{q^2}((1-q^2)q^{2M+4}s)].
$$
Due to Propositions \ref{p2.1} and \ref{p2.2} the last 
expression vanishes. 

If $s=1$ then 
$$ 
\int d_{q^2}ze_{q^2}(-i(1-q^2)z)E_{q^2}(i(1-q^2)q^2z)= 
$$ 
$$ 
=(1-q^2)\sum_{m=-\infty}^\infty q^{2m}
\left[\frac{(-i(1-q^2)q^{2m+2},q^2)_\infty}
{(-i(1-q^2)q^{2m},q^2)_\infty}+
\frac{(i(1-q^2)q^{2m+2},q^2)_\infty}
{(i(1-q^2)q^{2m},q^2)_\infty}\right]=
$$
$$
=(1-q^2)\sum_{m=-\infty}^\infty q^{2m}
\left(\frac1{1+i(1-q^2)q^{2m}}+\frac1{1-i(1-q^2)q^{2m}}\right)=
$$
$$
=2\sum_{m=\infty}^\infty\frac1{(1-q^2)^{-1}q^{-2m}+(1-q^2)q^{2m}}=
\frac2{1-q^2}\Theta_0
$$
(see (\ref{2.15}) and (\ref{2.17})).\rule{5pt}{5pt}

The next Lemma can be proved in the same way.
\begin{lem}\label{l5.2}
$$
\int\Ex E_{q^2}(i(1-q^2)q^2s)d_{q^2}s=
\left\{\begin{array}{lcl}
\frac2{1-q^2}\Theta_0 & for & z=1\\ 0 & for & z\ne 1.\\
\end{array}
\right.
$$
\end{lem}

\begin{lem}\label{l5.3} 
\begin{equation}\label{5.11} 
{\cal 
F}_{q^2}\Lambda=q^{-2}\Lambda^{-1}{\cal F}_{q^2}, 
\qquad {\cal F}_{q^2}\p_z=-is{\cal F}_{q^2}, 
\qquad {\cal F}_{q^2}z=-iq^{-2}\Lambda^{-1}\p_s{\cal F}_{q^2}, 
\end{equation} 
\begin{equation}\label{5.12}
{\cal F}_{q^2}^{-1}\Lambda=q^{-2}\Lambda^{-1}{\cal F}_{q^2}^{-1},
\qquad {\cal F}_{q^2}^{-1}\p_s=i\Lambda^{-1}z{\cal F}_{q^2}^{-1},
\qquad {\cal F}_{q^2}^{-1}s=i\p_z{\cal F}_{q^2}^{-1}.  
\end{equation} 
\end{lem}
{\bf Proof.} It is easy to verify that the kernels of the 
$q^2$-integral transforms (\ref{5.1}) and (\ref{5.4}) satisfy the
conditions:  
$$ 
\bg=\Eizso, 
$$ 
$$ 
\Ex=\eizso 
$$ 
and 
$$ 
\p_z\ddag E_{q^2}((1-q^2)azs)\ddag=a\ddag 
E_{q^2}((1-q^2)aq^2zs)\ddag s, 
$$ 
$$ 
\p_s\ddag E_{q^2}((1-q^2)azs)\ddag=az\ddag E_{q^2}((1-q^2)aq^2zs)\ddag,
$$
$$
\p_z\ddag e_{q^2}((1-q^2)azs)\ddag=a\ddag e_{q^2}((1-q^2)azs)\ddag s,
$$
$$
\p_s\ddag e_{q^2}((1-q^2)azs)\ddag=az\ddag e_{q^2}((1-q^2)azs)\ddag.
$$
From (\ref{4.6}) and the last relations we come to
the statement of Lemma.\rule{5pt}{5pt}

Now we can prove Proposition \ref{p5.1}.\\
From (\ref{5.11} we obtain
$$
{\cal F}_{q^2}z^k\p_z^l\phi(z)=(-i)^{k+l}q^{-2k}(\Lambda^{-1}\p_s)^ks^l
{\cal F}_{q^2}\phi(z).
$$
On the other hand
$$
\p_s^ks^l=(-1)^kq^{k(2l-k+1)}\sum_{j=0}^k(-1)^jq^{j(j-1)}
\frac
{(q^{-2l},q^2)}{(1-q^2)^{k-j}}
\left[\begin{array}{c}
k \\j
\end{array}
\right ]_{q^2}
s^{l-k+j}\p_s^j.
$$ 
Therefore, if $\phi(z)$ satisfies (\ref{4.1}) its image ${\cal F}_{q^2}\phi(z)$
satisfies (\ref{top}). It means, in particular, that the image of converged
sequence in  $\hat S_{q^2}$ converges in  $\hat S^{q^2}$. The statement 
concerning ${\cal F}_{q^2}^{-1}$ can be proved in the similar fashion using
(\ref{5.12}).

To prove the second part of Proposition we consider
 the action of the Fourier operators on the basis:
\begin{equation}\label{5.9} 
{\cal F}_{q^2}\circ{\cal F}_{q^2}^{-1}\hat\psi_n^\pm(s)=\hat\psi_n^\pm(s), 
\end{equation} 
\begin{equation}\label{5.10} 
{\cal F}_{q^2}^{-1}\circ{\cal F}_{q^2}\hat\phi_n^\pm(z)=\hat\phi_n^\pm(z).  
\end{equation} 
Consider the first relation 
$$ 
{\cal F}_{q^2}\circ{\cal F}_{q^2}^{-1}\hat\psi_n^+(s)= 
\frac1{2\Theta_0}\int d_{q^2}z
\Bigl(\int\Eizkm \hat\psi_n^+(\xi)d_{q^2}\xi\Bigr)\bg= 
$$ 
$$ 
=\frac{1-q^2}{2\Theta_0}q^{2n}\int d_{q^2}ze_{q^2}(-i(1-q^2)q^{2n}z)
\bg=
$$
$$
=\frac{1-q^2}{2\Theta_0}\int d_{q^2}ze_{q^2}(-i(1-q^2)z)
\phantom._0\Phi_1(-;0;q^2,i(1-q^2)q^{-2n+2}zs).
$$
It follows from Lemma \ref{l5.1} that
$$
{\cal F}_{q^2}\circ{\cal F}_{q^2}^{-1}\hat\psi_n^+(s)=
\left\{\begin{array}{lcl}
1 & for & s=q^{2n} \\ 0 & for & s\ne q^{2n}.\\
\end{array}
\right.
$$
and we come to (\ref{5.9}). In the similar way Lemma 
\ref{l5.2} leads to (\ref{5.10}).\rule{5pt}{5pt}

\section{The $q^2$-Fourier transform of $q^2$-distributions} 
\setcounter{equation}{0}

\begin{defi}\label{d6.1}
The $q^2$-Fourier transform of  $q^2$-distribution $f\in\hat S_{q^2}'$, 
is the $q^2$-distribution $g\in(\hat S^{q^2})'$ defined by the equality 
\begin{equation}\label{6.1}
<g,\psi>=<f,\phi>
\end{equation}
where $\phi(z)$ is an arbitrary function in $S_{q^2}$ and 
$\psi(s)\in\hat S^{q^2}$ is its 
$q^2$-Fourier transform.
\end{defi}
Suppose that the $q^2$-distribution $f$ corresponds to $f(z)$ and 
$zf(z)$ is absolutely $q^2$-integrable function. Let 
$\phi(z)= {\cal F}_{q^2}^{-1}\hat\psi(s)$. Then 
$$ 
<f,\phi>=\frac1{2\Theta_0}\int d_{q^2}z\bar f(z)\int\Ex\psi(s)d_{q^2}s= 
$$
$$
=\frac1{2\Theta_0}\int\overline{\int d_{q^2}zf(z)\mEx}
\psi(s)d_{q^2}s=<g,\psi>.
$$
It means that the $q^2$-distribution $g$ corresponds to the
function 
\begin{equation}\label{6.2} 
g(s)=\frac1{2\Theta_0}\int d_{q^2}zf(z)\mEx.  
\end{equation}

Similarly, if $g$ is defined by the absolutely $q^2$-integrable 
function $g(s)$ and $\hat\psi(s)={\cal F}_{q^2}\hat\phi(z)$, then 
$$ 
<g,\psi>=\int\int d_{q^2}z\phi(z)\bg\bar g(s)d_{q^2}s= 
$$ 
$$ 
=\int d_{q^2}z\phi(z)\overline{\int\bgm g(s)d_{q^2}s}=<f,\phi>, 
$$ 
i.e., $f$  corresponds to
\begin{equation}\label{6.3} 
f(z)=\int\bgm g(s)d_{q^2}s.  
\end{equation}

The $q^2$-Fourier transform of the $q^2$-distributions from 
$\hat S_{q^2}'$ we  denote by ${\cal F}_{q^2}'$. The next Proposition 
follows from (\ref{6.2}), (\ref{6.3}) and Lemma
\ref{l5.3} 
\begin{predl}\label{p6.1} 
There are the following commutation relations in the space of $q^2$ distributions:
\begin{equation}\label{6.4} 
{\cal F}_{q^2}'\Lambda=q^{-2}\Lambda^{-1}{\cal F}_{q^2}', \qquad {\cal 
F}_{q^2}'\p_z=-i\Lambda^{-1}s{\cal F}_{q^2}', \qquad {\cal 
F}_{q^2}'z=-i\p_s{\cal F}_{q^2}'.  
\end{equation} 
$$
({\cal F}_{q^2}')^{-1}\Lambda=q^{-2}\Lambda^{-1}({\cal F}_{q^2}')^{-1},
 ({\cal F}_{q^2}')^{-1}\p_s=iz({\cal F}_{q^2}')^{-1},
 ({\cal F}_{q^2}')^{-1}s=
iq^{-2}\Lambda^{-1}\p_z({\cal F}_{q^2}')^{-1}.
$$
\end{predl}

We define the $q^2$-Fourier transforms of some $q^2$-distributions. 

1. It follows from (\ref{3.8}) that
$$ 
{\cal F}_{q^2}'z^{-1}=\frac1{2\Theta_0}\int d_{q^2}zz^{-1}\mEx=
\frac i{\Theta_0}\Theta(s).
$$
Thus, 
\begin{equation}\label{6.6}
{\cal F}_{q^2}'z^{-1}=\hat i\sign s
=i(\theta_{q^2}^+-\theta_{q^2}^-)
\end{equation}
as the distribution over $\hat S^{q^2}$.

2. (\ref{6.4}), (\ref{6.6}) and (\ref{3.8}) give 
(\ref{4.13}) that
\begin{equation}\label{6.7}
{\cal F}_{q^2}'1={\cal F}_{q^2}'zz^{-1}=-i\p_s{\cal F}_{q^2}'z^{-1}=
\p_s(\theta_{q^2}^+-\theta_{q^2}^-)=2\delta_{q^2}.
\end{equation}

3.
$$
{\cal F}_{q^2}'(\theta_{q^2}^+-\theta_{q^2}^-)=
\frac1{2\Theta_0}\int_0^\infty d_{q^2}z\mEx-
\frac1{2\Theta_0}\int_{-\infty}^0d_{q^2}z\mEx=
$$
$$
=\frac{1-q^2}{2\Theta_0}\sum_{m=-\infty}^\infty q^{2m}
[e_{q^2}(i(1-q^2)q^{2m}s)-e_{q^2}(-i(1-q^2)q^{2m}s)]=
$$
$$
=\frac{i(1-q^2)}{\Theta_0}\sum_{m=-\infty}^\infty q^{2m}
\sin_{q^2}((1-q^2)q^{2m}s).
$$
Using (\ref{2.20}), we obtain
$$
{\cal F}_{q^2}'(\theta_{q^2}^+-\theta_{q^2}^-)=
\frac{is^{-1}}{\Theta_0}\lim_{M\to\infty}(1-q^2)s
\sum_{m=-M}^\infty q^{2m}\sin_{q^2}((1-q^2)q^{2m})=
$$
\begin{equation}\label{6.8}
=\frac{is^{-1}}{\Theta_0}\lim_{M\to\infty}[1-
\cos_{q^2}((1-q^2)q^{2M}]=\frac{is^{-1}}{\Theta_0}.
\end{equation}

4. It follows from (\ref{4.13}), (\ref{6.4}) and (\ref{6.8}) that
$$
{\cal F}_{q^2}'\delta=\frac12{\cal F}_{q^2}'\p_z
(\theta_{q^2}^+-\theta_{q^2}^-)=-\frac i2\Lambda^{-1}s{\cal F}_{q^2}'
(\theta_{q^2}^+-\theta_{q^2}^-)=\frac1{2\Theta_0}.
$$

5. From (\ref{4.10}) and (\ref{4.11}), (\ref{6.7}) and 
(\ref{6.8}) we obtain
$$
{\cal F}_{q^2}'\theta_{q^2}^+=\frac12{\cal F}_{q^2}' 
(\theta_{q^2}^+-\theta_{q^2}^-+1)=\frac{is^{-1}}{2\Theta_0}+\delta_{q^2}, 
$$
$$
{\cal F}_{q^2}'\theta_{q^2}^-=\frac12{\cal F}_{q^2}'
(-\theta_{q^2}^++\theta_{q^2}^-+1)=-\frac{is^{-1}}{2\Theta_0}+\delta_{q^2}. 
$$
\begin{predl}\label{p6.2}
For an arbitrary integer $n\ge0$
$$
{\cal F}_{q^2}'z^n=2i^nq^{-n(n+1)}\frac{(q^2,q^2)_n}
{(1-q^2)^n}s^{-n}\delta_{q^2}(s).
$$
\end{predl}
{\bf Proof.} The relation 
$1=\frac{(1-q^2)^n}{(q^2,q^2)_n}\p_z^nz^n$ follows from (\ref{3.2}). 
Next, from 
(\ref{6.7}) and Proposition \ref{p6.1} we obtain
$$ 
\delta_{q^2}(s)=\frac12{\cal F}_{q^2}'1=\frac{(1-q^2)^n}
{2(q^2,q^2)_n} {\cal F}_{q^2}'\p_z^nz^n=(-i)^n\frac{(1-q^2)^n}
{2(q^2,q^2)_n}(\Lambda^{-1}s)^n{\cal F}_{q^2}'z^n.
$$
Thus,
\begin{equation}\label{6.13}
{\cal F}_{q^2}'z^n=2i^n\frac{(q^2,q^2)_n}{(1-q^2)^n}
(s^{-1}\Lambda)^n\delta_{q^2}(s).
\end{equation}
Since $(s^{-1}\Lambda)^n=q^{-n(n-1)}s^{-n}\Lambda^n$ and 
$\delta_{q^2}(q^{2n}s)=q^{-2n}\delta_{q^2}(s)$, 
$$
(s^{-1}\Lambda)^n\delta_{q^2}(s)=q^{-n(n+1)}s^{-n}\delta_{q^2}(s).
$$
Then, the statement of Proposition follows from  (\ref{6.13}). 
\rule{5pt}{5pt}

\begin{predl}\label{p6.3}
For an arbitrary integer $n\ge0$
$$
{\cal F}_{q^2}'z^{-n-1}=
i^{n+1}\frac{(q^2,q^2)_n}{(1-q^2)^n}s^n\sign s.
$$
\end{predl}
{\bf Proof.} For an arbitrary $n\ge0$ we have from (\ref{3.3}) 
$$
\p_z^nz^{-1}=(-1)^nq^{-n(n+1)}\frac{(q^2,q^2)_n}{(1-q^2)^n}z^{-n-1}.
$$
Using Proposition \ref{p6.1} we obtain
$$
{\cal F}_{q^2}'z^{-n-1}=(-1)^nq^{n(n+1)}\frac{(1-q^2)^n}{(1-q^2)^n}
{\cal F}_{q^2}'\p_z^nz^{-1}=
i^nq^{n(n+1)}\frac{(1-q^2)^n}{(1-q^2)^n}(\Lambda^{-1}s)^n
{\cal F}_{q^2}'z^{-1}.
$$
It can be derived from (\ref{4.6}) that 
$(\Lambda^{-1}s)^n=q^{-n(n+1)}s^n\Lambda^{-n}$.
Using this equality and (\ref{6.6}), we come to the statement of 
Proposition. \rule{5pt}{5pt} 
\begin{predl}\label{p6.4}
\begin{equation}\label{6.15}
{\cal F}_{q^2}'z_+^{\nu-1}=
\frac{e_{q^2}(q^2)E_{q^2}(-q^{2(1-\nu)})}{2\Theta_0}
(\bar c_\nu s_-^{-\nu}+c_\nu s_+^{-\nu}),
\end{equation}
\begin{equation}\label{6.16}
{\cal F}_{q^2}'z_-^{\nu-1}=
-\frac{e_{q^2}(q^2)E_{q^2}(-q^{2(1-\nu)})}{2\Theta_0}
(c_\nu s_-^{-\nu}+\bar c_\nu s_+^{-\nu}),
\end{equation}
where
$$
c_\nu=\sum_{m=-\infty}^\infty\frac{q^{2\nu m}(q^{-2m}+i(1-q^2))}
{(1-q^2)^{-1}q^{-2m}+(1-q^2)q^{2m}}.
$$
\end{predl}
{\bf Proof.} If $0<\R\nu<1$,  (\ref{4.15}) and (\ref{6.2}) give
\begin{equation}\label{6.18}
{\cal F}_{q^2}z_+^{\nu-1}=
\frac{1-q^2}{2\Theta_0}\sum_{m=-\infty}^\infty q^{2\nu m}
e_{q^2}(i(1-q^2)q^{2m}s).
\end{equation}
It can be directly checked that $z_+^{\nu-1}$ satisfies the equation
$$
z\p_zf(z)=\frac{1-q^{2(\nu-1)}}{1-q^2}f(z).
$$
From the Proposition \ref{p6.1} we obtain that  
the $q^2$-Fourier transformed equation takes the form
$$ 
-q^{-2}\p_s(sg(q^{-2}s))=\frac{1-q^{2(\nu-1)}}{1-q^2}g(s).  
$$ 
Since $s_-^{-\nu}$ and $s_+^{-\nu}$ satisfy this equation we can write
\begin{equation}\label{6.19}
{\cal F}_{q^2}'z_+^{\nu-1}=c_1s_-^{-\nu}+c_2s_+^{-\nu}.
\end{equation}
Substituting $s=1$ in (\ref{6.19}) we obtain from (\ref{2.8}),(\ref{6.18}) and 
(\ref{6.19}) 
$$ 
c_2=\frac{1-q^2}{2\Theta_0}\sum_{m=-\infty}^\infty q^{2\nu 
m}e_{q^2}(i(1-q^2)q^{2m})=
$$
$$
=\frac{1-q^2}{2\Theta_0}e_{q^2}(q^2)\sum_{m=-\infty}^\infty q^{2\nu m}
\sum_{n=0}^\infty\frac{(-1)^nq^{n(n+1)}}
{(q^2,q^2)_n[1-i(1-q^2)q^{2(m+n)}]}.
$$
If $0<\R\nu<1$, we can change the order of summations. Then
$$ 
c_2=\frac{1-q^2}{2\Theta_0}e_{q^2}(q^2)\sum_{n=0}^\infty 
\frac{(-1)^nq^{n(n+1)}}{(q^2,q^2)_n}\sum_{m=-\infty}^\infty
\frac{q^{2\nu m}}{1-i(1-q^2)q^{2(m+n)}}=
$$
$$
=\frac{e_{q^2}(q^2)}{2\Theta_0}\sum_{n=0}^\infty
\frac{(-1)^nq^{n(n-1)}q^{2(1-\nu)n}}{(q^2,q^2)_n}
\sum_{m=-\infty}^\infty\frac{q^{2\nu m}(q^{-2m}+i(1-q^2))}
{(1-q^2)^{-1}q^{-2m}+(1-q^2)q^{2m}}.
$$
We denote the sum of the inner series by $c_\nu$. Then
$$
c_2=\frac{c_\nu}{2\Theta_0}e_{q^2}(q^2)E_{q^2}(-q^{2(1-\nu)})
$$
Assume now that $s=-1$ in (\ref{6.19}). Then in the same way we obtain
$$ 
c_1=\frac{\bar c_\nu}{2\Theta_0}e_{q^2}(q^2)E_{q^2}(-q^{2(1-\nu)}) 
$$ 
and thereby come to
(\ref{6.15}). It determines ${\cal F}_{q^2}'z_+^{\nu-1}$ as  
(\ref{6.15}) for non-integer $\nu$ by the analytical continuation.

The $q^2$-Fourier transform of $z_-^{\nu-1},\quad \nu\ne0,\pm1,\ldots.$
$$
{\cal F}_{q^2}'z_-^{\nu-1}=
-\frac{1-q^2}{2\Theta_0}\sum_{m=-\infty}^\infty q^{2\nu m}
e^{q^2}(-i(1-q^2)q^{2m}s).
$$
 differs from (\ref{6.18}) 
 by the sign and the complex conjugation, i.e., in this way
 we have (\ref{6.16}).
\rule{5pt}{5pt}

\bigskip

 We summarize  all results in the table:

\bigskip
\begin{tabular}{|p{3cm}|p{12cm}|}
\hline
$f(z)$ & $g(s)=\frac1{2\Theta_0}\int d_{q^2}zf(z)\mEx$\\ 
\hline
\hline 
$\delta_{q^2}(z)$ & $\frac1{2\Theta_0}$\\ 
\hline 
$\theta_{q^2}^+(z)$ & $\frac i{2\Theta_0}s^{-1}+\delta_{q^2}(s)$\\ 
\hline 
$\theta_{q^2}^-(z)$ & $-\frac i{2\Theta_0}s^{-1}+\delta_{q^2}(s)$\\ 
\hline 
$z^n$ & $2i^nq^{-n(n+1)}\frac{(q^2,q^2)_n}
{(1-q^2)^n}s^{-n}\delta_{q^2}(s)$\\
$(n=0,1,\ldots)$ & \\
\hline
$z^{-n-1}$ & $i^{n+1}\frac{(1-q^2)^n}{(q^2,q^2)_n}s^n\sign s$\\
$(n=0,1,\ldots)$ & \\
\hline
$z_+^{\nu-1}$ & $\frac{e_{q^2}(q^2)E{q^2}(-q^{2(1-\nu)})}{2\Theta_0}
(\bar c_\nu s_-^{-\nu}+c_\nu s_+^{-\nu})$\\
$(\nu\ne0,\pm1,\ldots)$ & \\
\hline
$(z_-^{\nu-1}$ &
$-\frac{e_{q^2}(q^2)E{q^2}(-q^{2(1-\nu)})}{2\Theta_0}
(c_\nu s_-^{-\nu}+\bar c_\nu s_+^{-\nu}]$\\
$(\nu\ne0,\pm1,\ldots)$ & \\
\hline
\end{tabular}

{\bf Acknowledgments.}\\
{\sl  
 The work  is 
supported in part by grants 96-15-96455 for support of scientific schools,
INTAS-93-0166 extension and by Award No.
RM2-150 of the Civilian Research \& Development Foundation
(CRDF) for the Independent States of the Former Soviet Union
(M.O.); by RFBR 97-01-00747 and NIOKR MPS RF  (V.R.)}
The work is prepared partly during visit of M.O. in IHES (Bur-sur-Yvette). 
He is grateful to Prof. K.Gawedzky for hospitality.

\bigskip
\small{

}

\end{document}